\documentclass[12pt]{article}
\usepackage{amssymb} 
\usepackage{amsmath} 
\title{Potential modularity---a survey}
\author{Kevin Buzzard}
\newcommand{\Gal}{\mathop{\mathrm{Gal}}}
\newcommand{\PGL}{\mathop{\mathrm{PGL}}}
\newcommand{\re}{\mathop{\mathrm{Re}}}
\newcommand{\GL}{\mathop{\mathrm{GL}}}
\newcommand{\univ}{{\mathop{\mathrm{univ}}}}
\newcommand{\pst}{{\mathop{\mathrm{pst}}}}
\newcommand{\Frac}{{\mathop{\mathrm{Frac}}}}
\newcommand{\Spec}{{\mathop{\mathrm{Spec}}}}
\newcommand{\calO}{\mathcal{O}}
\newcommand{\m}{\mathfrak{m}}
\newcommand{\Qbar}{\overline{\Q}}
\newcommand{\Q}{\mathbf{Q}}
\newcommand{\Qp}{\Q_p}
\newcommand{\Qpbar}{\Qbar_p}
\newcommand{\rhobar}{\overline{\rho}}
\newcommand{\R}{\mathbf{R}}
\newcommand{\GQ}{\Gal(\Qbar/\Q)}
\newcommand{\GQp}{\Gal(\Qpbar/\Qp)}
\newcommand{\GF}{\Gal(\overline{F}/F)}
\newcommand{\F}{\mathbf{F}}
\newcommand{\Fpbar}{\overline{\F}_p}
\newcommand{\Flbar}{\overline{\F}_\ell}
\newcommand{\p}{\mathfrak{p}}
\newcommand{\T}{\mathbf{T}}
\newcommand{\Z}{\mathbf{Z}}
\usepackage{amsthm}
\theoremstyle{plain}
\newtheorem{theorem}{Theorem}[section]
\newtheorem{lemma}[theorem]{Lemma}

\newtheorem{conjecture}[theorem]{Conjecture}

\theoremstyle{remark}

\begin{document}
\maketitle
\begin{abstract}
A Spitalfields Day at the Newton Institute was organised on the
subject of the recent theorem that any elliptic curve over any
totally real field is potentially modular. This article is a
survey of the strategy of the proof, together with some history.
\end{abstract}

\section*{Introduction}

Our main goal in this article is to talk about recent theorems
of Taylor and his co-workers on modularity and potential modularity
of Galois representations, particularly those attached to elliptic curves.
However, so as to not bog down the exposition unnecessarily with technical
definitions right from the off, we will build up to these results by starting
our story with Wiles' breakthrough paper~\cite{wiles:flt}, and working towards
the more recent results. We will however assume some familiarity with the
general area---for example we will assume the reader is familiar with the
notion of an elliptic curve over a number field, and a Galois representation,
and what it means for such things to be modular (when such a notion
makes sense). Let us stress now that, because of this chronological
approach, some theorems stated in this paper will be superseded
by others (for example Theorem~\ref{wtw} gets superseded by Theorem~\ref{cdtt}
which gets superseded by Theorem~\ref{stw}), and similarly some conjectures
(for example Serre's conjecture) will become theorems as the story progresses.
The author hopes that this slightly non-standard style nevertheless
gives the reader the feeling of seeing how the theory evolved. 

We thank Toby Gee for reading through a preliminary draft of this article
and making several helpful comments, and we also thank Matthew Emerton
and Jan Nekov{\'a}{\v{r}} for pointing out various other inaccuracies
and ambiguities.

\section{Semistable elliptic curves over $\Q$ are modular.}

The story, of course, starts with the following well-known result
proved in~\cite{wiles:flt} and~\cite{tw:flt}.
\begin{theorem}[Wiles, Taylor--Wiles]\label{wtw} Any semi-stable elliptic curve
over the rationals is modular.
\end{theorem}
This result, together with work of Ribet and others on Serre's conjecture,
implies Fermat's Last Theorem. This meant that the work of Wiles and Taylor
captured the imagination of the public. But this article is not about Fermat's
Last Theorem, it is about how the \emph{modularity theorem} above
has been vastly generalised. Perhaps we should note here though that
that there is still a long way to go! For example, at the time of
writing, it is still an open
problem as to whether an arbitrary elliptic curve over an arbitrary
totally real field is modular, and over an general number field,
where we cannot fall back on the theory of Hilbert modular forms,
the situation is even worse (we still do not have a satisfactory
theorem attaching elliptic curves to modular forms in this generality,
let alone a result in the other direction).

Before we go on to explain the generalisations that this article
is mainly concerned with, we take some time to remind the reader
of some of the details of the strategy of the Wiles/Taylor--Wiles proof.
The ingredients are as follows. For the first, main, ingredient,
we need to make some definitions. Let $p$ be a prime number,
let $\calO$ denote the integers in a finite extension of $\mathbf{Q}_p$,
let $\rho:\Gal(\Qbar/\Q)\to\GL_2(\calO)$ denote a continuous odd
(by which we mean $\det(\rho(c))=-1$, for $c$ a complex conjugation)
irreducible representation, unramified outside a finite set
of primes, and let $\rhobar$
denote its reduction modulo the maximal ideal of~$\calO$. Recall
that there is a general theorem due to Deligne and others,
which attaches $p$-adic Galois representations to modular eigenforms;
we say a $p$-adic Galois representation~$\rho$ is \emph{modular}
if it arises in this way, and that a mod~$p$ Galois representation~$\rhobar$
is \emph{modular} if it arises as the semisimplification
of the mod $p$ reduction of a modular $p$-adic Galois representation. 
Note that we will allow myself the standard abuse of notation here,
and talk about ``mod $p$ reduction'' when we really mean ``reduction
modulo the maximal ideal of $\calO$''.

Note that if a $p$-adic Galois representation $\rho$ is modular, then its
reduction $\rhobar$
(semisimplified if necessary) is trivially modular. Wiles'
insight is that one could sometimes go the other way.

\begin{theorem}[Modularity lifting theorem]\label{mlt1} If $p>2$, if $\rhobar$
is irreducible and modular, and if furthermore $\rho$ is semistable and
has cyclotomic determinant, then $\rho$ is modular.
\end{theorem}
This is Corollary~3.46 of~\cite{ddt} (the aforementioned paper
is an overview of the Wiles/Taylor--Wiles work; the theorem
is essentially due to Wiles and Taylor--Wiles). The proof is some hard
work, but is now regarded as ``standard''---many mathematicians
have read and verified the proof. We shall say a few words about
the proof later on. Semistability
is a slightly technical condition (see op.\ cit.\ for
more details) but we shall be removing it soon so we do
not go into details. Rest assured that if $E/\Q$ is a semistable
elliptic curve then its Tate module is semistable.

The next ingredient is a very special case of the following conjecture
of Serre:

\begin{conjecture}[Serre, 1987]\label{serreconj} If $k$ is a finite
field and $\rhobar:\GQ\to\GL_2(k)$ is a continuous odd absolutely
irreducible representation, then $\rhobar$ is modular.
\end{conjecture}

We will say more about this conjecture and its generalisations later.
Note that this conjecture is now a theorem of Khare and Wintenberger,
but we are taking a chronological approach so will leave it as
a conjecture for now. In the early 1990s this conjecture was wide open, but one
special case had been proved in 1981 (although perhaps it was not
stated in this form in 1981; see for example Prop.~11 of~\cite{serre:conj}
for the statement we need):

\begin{theorem}[Langlands, Tunnell]\label{lt} If $\rhobar:\GQ\to\GL_2(\F_3)$
is continuous, odd, and irreducible, then $\rhobar$ is modular.
\end{theorem}
The original proof of Theorem~\ref{lt}
is a huge amount of delicate analysis: let it not
be underestimated! One needs (amongst other things)
the full force of the trace formula
in a non-compact case to prove this result, and hence a lot of delicate
analysis. We freely confess to not having
checked the details of this proof ourselves. Note also that the point is
\emph{not} just that the image of $\rhobar$ is solvable, it is that the image
is very small (just small enough to be manageable, in fact).
Note that because we are in odd characteristic, the notions of
irreducibility and absolute irreducibility coincide for an odd
representation (complex conjugation has two distinct eigenvalues,
both defined over the ground field).
Langlands' book~\cite{langlands:bc} proves
much of what is needed; the proof was finished by Tunnell
in~\cite{tunnell:artin} using the non-solvable cubic base
change results of~\cite{jpss}. This result is of course also
now regarded as standard---many mathematicians have read
and verified this proof too. The author remarks however that due to the rather
different techniques involved in the proofs of the two results,
he has the impression that the number of mathematicians who have
read and verified all the details of
the proofs of \emph{both} the preceding theorems
is rather smaller!

Let us see how much of Theorem~\ref{wtw} we can prove so far, given
Theorems~\ref{mlt1} and~\ref{lt}.
Let $E$ be a semistable elliptic curve, set $p=3$ and
let $\rho:\GQ\to\GL_2(\Z_3)$ be the 3-adic Tate module of $E$.
Then $\rho$ is continuous, odd, unramified outside a finite
set of primes, and semistable, with cyclotomic determinant.
Furthermore, \emph{if} $\rhobar$, the Galois representation on $E[3]$,
is irreducible, then $\rhobar$ is modular by the Langlands--Tunnell
theorem~\ref{lt} and
so $\rho$, and hence~$E$, is modular, by the modularity lifting
theorem. Of course the problem is that $E[3]$ may be reducible---for
example if~$E$ has a $\Q$-point of order~3 (or more generally a subgroup
of order~3 defined over~$\Q$). To deal with this situation, Wiles
developed a technique known as the ``3--5 trick''.

\begin{lemma}[The 3--5 trick]\label{35trick} If $E/\Q$ is a semistable elliptic curve
with $E[3]$ reducible, then $E[5]$ is irreducible, and
there is another semistable elliptic curve $A/\Q$
with $E[5]\cong A[5]$ and $A[3]$ irreducible.
\end{lemma}

The 3--5 trick is all we need to finish the proof of Theorem~\ref{wtw}.
For if $E/\Q$ is semistable but $E[3]$ is reducible, choose~$A$ as
in the lemma, and note that $A[3]$ is irreducible, hence $A[3]$
is modular (by the Langlands--Tunnell Theorem~\ref{lt}), hence $A$ is
modular (by the Modularity Lifting Theorem~\ref{mlt1}), so $A[5]$
is modular, so $E[5]$
is modular, and irreducible, so $E$ is modular by the Modularity
Lifting Theorem~\ref{mlt1} applied to the 5-adic Tate module of~$E$.

Let us say a few words about the proof of the last lemma. The reason
that reducibility of $E[3]$ implies irreducibility of $E[5]$
is that reducibility of both would imply that $E$ had a rational
subgroup of order~15, and would hence give rise to a point
on the modular curve $Y_0(15)$, whose compactification $X_0(15)$
is a curve of genus~1 with
finitely many rational points, and it turns out that the points
on this curve are known, and one can check that none of them
can come from semistable elliptic curves (and all of them
are modular anyway). So what is left is that given~$E$ as
in the lemma, we need to produce~$A$. Again we use
a moduli space trick. We consider the moduli space over~$\Q$
parametrising elliptic
curves~$B$ equipped with an isomorphism $B[5]\cong E[5]$ that preserves
the Weil pairing. This moduli space has a natural smooth compactification~$X$
over~$\Q$, obtained by adding cusps. Over the complexes the resulting
compactified curve is isomorphic to the modular curve $X(5)$, which
has genus zero. Hence~$X$ is a genus zero curve over~$\Q$. Moreover,
$X$ has a rational point (coming from~$E$) and hence $X$ is itself
isomorphic to the projective line over~$\Q$ (rather than a twist
of the projective line). In particular, $X$ has infinitely many
rational points. Now using Hilbert's Irreducibility Theorem,
which in this setting can be viewed as some sort of refinement
of the Chinese Remainder Theorem, it is possible
to find a point on~$X$ which is 5-adically very close to~$E$,
and 3-adically very far away from~$E$ (far enough so that
the Galois representation on the 3-torsion of the corresponding
elliptic curve is irreducible: this is the crux of the Hilbert
Irreducibility Theorem, and this is where we are using more
than the naive Chinese Remainder Theorem). Such a point
corresponds to the elliptic curve~$A$ we seek, and the lemma,
and hence Theorem~\ref{wtw}, is proved.

The reason we have broken up the proof of Theorem~\ref{wtw} into
these pieces is that we would like to discuss generalisations
of Theorem~\ref{wtw}, and this will entail discussing generalisations
of the pieces that we have broken it into.

\section{Why the semistability assumption?}

All semistable elliptic curves were known to be modular
by 1995, but of course one very natural question was whether
the results could be extended to all elliptic curves. Let us
try and highlight the issues involved with trying to extend
the proof; we will do this by briefly reminding readers of the strategy
of the \emph{proof} of a modularity lifting theorem such
as Theorem~\ref{mlt1}. The strategy is
that given an irreducible modular $\rhobar$, one considers two kinds of lifting
to characteristic zero. The first is a ``universal deformation''
$\rho^{\univ}:\GQ\to\GL_2(R^{\univ})$, where one considers \emph{all} deformations
satisfying certain properties fixed beforehand
(in Wiles' case these properties were typically
``unramified outside $S$ and semistable at all the primes in $S$'' for some
fixed finite set of primes~$S$), and uses the result of Mazur
in~\cite{mazur:deformations} that says that there
is a \emph{universal} such deformation, taking values in a ring $R^{\univ}$.
The second is a ``universal modular deformation''
$\rho_{\T}:\GQ\to\GL_2(\T)$ comprising of a lift of $\rhobar$ to
a representation taking values in a Hecke algebra over $\Z_p$ built from
modular forms of a certain level, weight and character (or perhaps
satisfying some more refined local properties). Theorems about
modular forms (typically local-global theorems) tell us that the
deformation to $\GL_2(\T)$ has the properties used in the definition of
$R^{\univ}$, and there is hence a map
$$R^{\univ}\to\T.$$
The game is to prove that this map is an isomorphism; then
all deformations will be modular, and in particular $\rho$, the
representation we started with, will be modular. The insight
that the map may be an isomorphism seems to be due to Mazur: see
Conjecture (*) of~\cite{mazur-tilouine} and the comments preceding it.
One underlying miracle is
that this procedure can only work if $R^{\univ}$ has no $p$-torsion,
something which is not at all evident, but which came out
of the Wiles/Taylor--Wiles proof as a consequence.

The original proof that the map $R^{\univ}\to\T$ is an isomorphism breaks
up into two steps: the first one, referred to as the minimal case,
deals with situations where $\rho$ is ``no more ramified than
$\rhobar$'', and is proved by a patching argument via the
construction of what is now known as a Taylor--Wiles system: one
checks that certain projective limits of $R$s and $T$s (using weaker
and weaker deformation conditions, and more and more modular forms)
are power series rings, and that the natural map between them is an isomorphism
for commutative algebra reasons (for dimension reasons, really),
and then one descends back to the case of interest. The second
is how to deduce the general case from the minimal case---this
is an inductive procedure (which relies on a result on Jacobians
of modular curves known as Ihara's Lemma, the analogue
of which still appears to be open for $\GL_n$, $n>2$; this provided
a serious stumbling-block in generalising the theory to higher dimensions
for many years). Details of both of these arguments can be found
in~\cite{ddt}, especially \S5 (as well, of course, as in the original
sources).

To see why the case of semistable elliptic curves was treated
first historically, we need to look more closely at the nitty-gritty
of the details behind a deformation problem.

Wiles had a semistable irreducible representation $\rhobar:\GQ\to\GL_2(k)$ with
$k$ a finite field of characteristic~$p$.
Say $\rhobar$ is unramified outside some finite
set~$S\ni p$ of primes, and has cyclotomic determinant. Crucially, Wiles
knew what it meant (at least when $p>2$) for a deformation $\rho:\GQ\to\GL_2(A)$
to be semistable and unramified outside~$S$, where $A$ is now a general
Artin local ring with residue field~$k$ (or even a projective limit
of such rings). For a prime $q\not\in S$ it of course means $\rho$
is unramified at~$q$. For $q\not=p$, $q\in S$, it means
that the image of an inertia group at~$q$ under $\rho$ can be conjugated
into the upper triangular unipotent matrices. For $q=p$ one needs
more theory. The observation is that an elliptic curve with
semistable reduction either has good reduction, or multiplicative
reduction. The crucial point is that for a general Artin local~$A$
with finite residue field one can make sense of the
notion that $\rho:\GQ\to\GL_2(A)$ has ``good reduction''---one demands that it
is the Galois action on the generic fibre of a finite flat group
scheme with good reduction, and work of Fontaine~\cite{fontaine:book}
and Fontaine--Laffaille~\cite{fl} shows that one can translate this notion
into ``linear algebra'' which
is much easier to work with (this is where the assumption $p>2$ is needed).
Also crucial are the results
of Raynaud~\cite{ppp}, which show that the category of Galois
representations with these properties is very well-behaved.
Similarly one can make sense of the notion that $\rho$ has
``multiplicative reduction'': one can demand that $\rho$
on a decomposition group at~$p$ is upper triangular.

We stress again that the crucial point is that the notions of ``good reduction''
and ``multiplicative reduction'' above make sense for an \emph{arbitrary}
Artin local~$A$, and patch together well to give well-behaved local
deformation conditions which are locally representable (by which we mean
the deformations of the Galois representation $\rhobar|G_{\Qp}$
are represented by some universal ring). So we get
a nicely-behaved local deformation ring---in particular we get
a ring for which we can compute the tangent space $\m/\m^2$
of its mod~$p$ reduction. If this tangent space has dimension
at most~1 then the dimension calculations work out in the patching
argument and the modularity lifting theorem follows.

If one is prepared to take these observations on board, then it
becomes manifestly clear what the one of the main problems will be in proving
that an arbitrary elliptic curve over the rationals is modular:
we will have to come up with deformation conditions that are small
enough to make the dimension calculations work, but big enough to encompass
Galois representations that are not semistable. At primes $q\not=p$
this turned out to be an accessible problem; careful calculations
by Fred Diamond in~\cite{fred:15} basically resolved these issues
completely. Diamond's main theorem had as a consequence the result
that if~$E/\Q$ had semistable
reduction at both~3 and~5 then~$E$ was modular. The reason that
both~3 and~5 occur is of course because he has to use the 3--5 trick
if $E[3]$ is reducible. A year or so later, Diamond, and independently
Fujiwara, had another insight: instead of taking limits of Hecke
algebras to prove that a deformation ring equalled a Hecke ring,
one can instead take limits of modules that these algebras act on
naturally. The resulting commutative algebra is more delicate,
and one does not get modularity of any more elliptic curves
in this way, but the result is of importance because it enables
one to apply the machinery in situations where certain ``mod~$p$
multiplicity one'' hypotheses are not known. These multiplicity one
hypotheses were known in the situations that Wiles initially dealt with
but were not known in certain more general situations; the consequence
was that Wiles' method could now be applied more generally. 
Diamond's paper~\cite{diamond:mult1} illustrated the point by showing
that the methods could now be applied in the case of Shimura
curves over~$\Q$ (where new multiplicity one results could be deduced
as a byproduct), and Fujiwara (in~\cite{fuj:mult1}, an article
which remains unpublished, for reasons unknown to this author)
illustrated that the method enabled one to generalise Wiles'
methods to the Hilbert modular case, on which more later.

Getting back to elliptic curves over the rationals, the situation
in the late 1990s, as we just indicated, was that any elliptic
curve with semistable reduction at~3 and~5 was now proven
to be modular. To get further, new ideas were needed, because in the 1990s
the only source of modular mod~$p$ Galois representations were
those induced from a character, and those coming from the
Langlands--Tunnell theorem. Hence in the 1990s one was forced to
ultimately work with the
prime $p=3$ (the prime $p=2$ was another possibility; see for
example~\cite{dickinson}, but here other technical issues
arise). Hence, even with the 3--5 trick, it was clear that if one
wanted to prove that all elliptic curves over $\Q$ were modular
using these methods then one was going to have
to deal with elliptic curves that have rather nasty non-semistable
reduction at~3 (one cannot use the 3--5 trick to get around this
because if $E$ has very bad reduction at~3 (e.g.\ if its conductor
is divisible by a large power of~3) then this will be reflected
in the 5-torsion, which will also have a large power of~3 in its
conductor, so any curve~$A$ with $A[5]\cong E[5]$ will also be
badly behaved at~3; one can make certain simplifications this way
but one cannot remove the problem entirely). The main
problem is then deformation-theoretic: given some elliptic
curve~$E$ which is highly ramified at some odd prime~$p$, how does
one write down a reasonable deformation problem for~$E[p]$ at~$p$, which
is big enough to see the Tate module of the curve, but is still sufficiently
small for the Taylor--Wiles method to work? By this we mean that the
tangent space of the mod~$p$ local deformation problem at~$p$ has to have
dimension at most~1. This thorny issue explains the five year gap
between the proof of the modularity of all semistable elliptic curves,
and the proof for all elliptic curves.

\section{All elliptic curves over~$\Q$ are modular.}

As explained in the previous section, one of the main obstacles in
proving that all elliptic curves over the rationals are modular
is that we are forced, by Langlands--Tunnell, to work with
$p=3$, so elliptic curves with conductor a multiple of a high power
of~3 are going to be difficult to deal with. Let us review the situation
at hand. Let $\rhobar:\GQp\to\GL_2(k)$ be an irreducible
Galois representation, where here $k$ is a finite field
of characteristic~$p$. Such a $\rhobar$ has a universal
deformation to $\rho^\univ:\GQp\to\GL_2(R^\univ)$. This
ring $R^\univ$ is a quotient of a power series ring
$W(k)[[x_1,x_2,\ldots,x_n]]$ in finitely many variables,
where here $W(k)$ denotes the Witt vectors of the field~$k$.
But this universal deformation ring is too big for our purposes---a general
lifting of $\rhobar$ to the integers $\calO$ of a finite extension
of $\Frac(W(k))$ will not look anything like the Tate module
of an elliptic curve (it will probably not even be Hodge--Tate,
for example). The trick that Wiles used was to not look at
such a big ring as $R^{\univ}$, but to look at more stringent
deformation problems, such as deforming $\rhobar$ to representations
which came from finite flat group schemes over $\Z_p$. This
more restricted space of deformations is represented
by a smaller deformation ring $R^\flat$, a quotient of $R^\univ$,
and it is rings such as $R^\flat$ that Wiles could work with
(the relevant computations in this case were done in Ravi
Ramakrishna's thesis~\cite{ramakrishna:thesis}).

In trying to generalise this idea we run into a fundamental problem.
The kind of deformation problems that one might want to look
at are problems of the form ``$\rho$ that become finite and
flat when restricted to $\Gal(\Qpbar/K)$ for $K$ this fixed
finite extension of $\Q_p$''. However, for $K$ a wildly ramified
extension of $\Q_p$ the linear algebra methods alluded to earlier
on become much more complex, and indeed at this point historically
there was no theorem classifying finite flat group schemes
over the integers of such $p$-adic fields which was concrete
enough to enable people to check that the resulting deformation
problems were representable, and represented by rings whose tangent spaces were
sufficiently small enough to enable the methods to work.

A great new idea, however, was introduced in the paper~\cite{cdt}. 
Instead of trying to write down a complicated deformation problem
that made sense for all Artin local rings and then to analyse the
resulting representing ring, Conrad, Diamond and Taylor construct
``deformation rings'' in the following manner. First, they consider
deformations $\rho:\GQp\to\GL_2(\calO)$ of $\rhobar$ in the
case that $\calO$ is the integers of a finite extension of $\Q_p$.
In this special setting there is a lot of extra theory available: one
can ask if $\rho$ is Hodge--Tate, de Rham, potentially semi-stable,
crystalline and so on (these words \emph{do not make sense} when
applied to a general deformation $\rho:\GQp\to\GL_2(A)$
of $\rhobar$, they only make
sense when applied to a deformation to $\GL_2(\calO)$),
and furthermore if $\rho$ is potentially
semi-stable then the associated Fontaine module $D_\pst(\rho)$
is a 2-dimensional vector space with an action of the inertia
subgroup of $\GQp$ which
factors through a finite quotient. This finite image
2-dimensional representation of inertia is called the \emph{type} of
the potentially semi-stable representation $\rho$, and so we
can fix a type $\tau$ and then ask
that a deformation $\rho:\GQp\to\GL_2(\calO)$ be potentially semistable of
a given type.

Again we stress that this notion of being potentially semistable
of a given type certainly
does \emph{not} make sense for a deformation of $\rhobar$
to an arbitrary Artin local $W(k)$-algebra, so in particular
this notion is \emph{not} a deformation problem and we cannot
speak of its representability. One of the insights of~\cite{cdt}
however, is that we can construct a ``universal ring''
for this problem anyway! Here is the trick, which is really
rather simple. We have $\rhobar$ and its universal formal
deformation $\rho^\univ$ to $GL_2(R^\univ)$. Now let us consider
all maps $s:R^\univ\to\calO$, where $\calO$ is as above.
Given such a map $s$, we can compose $\rho^\univ$ with $s$
to get a map $\rho_s:\GQp\to\GL_2(\calO)$. Let us say that the kernel
of $s$ is \emph{of type~$\tau$} if $\rho_s$ is of type~$\tau$,
and if furthermore $\rho_s$ is potentially Barsotti--Tate (that
is, comes from a $p$-divisible group over the integers of a finite
extension of $\Q_p$) and has determinant
equal to the cyclotomic character. A good example of a potentially Barsotti--Tate
representation is the representation coming from the Tate
module of an elliptic curve with potentially good reduction at~$p$,
and such things will give rise to points of type~$\tau$ for an appropriate
choice of~$\tau$.

Let $R_\tau$ denote the quotient of $R^\univ$ by the intersection
of all the prime ideals of $R^\univ$ which are of type $\tau$
(with the convention that $R_\tau=0$ if there are no such prime ideals).
Geometrically, what is happening is that the kernel of $s$
is a prime ideal and hence a point in $\Spec(R^\univ)$, and
we are considering the closed subscheme of $\Spec(R^\univ)$
obtained as the Zariski-closure of all the prime ideals of type $\tau$.
So, whilst $R_\tau$ does not represent the moduli problem of
being ``of type $\tau$'' (because this is not even a moduli problem,
as mentioned above), it is a very natural candidate for a ring
to look at if one wants to consider deformations of type $\tau$.
It also raises the question as to whether the set of points
which are of type $\tau$ actually form a closed set in,
say, the rigid space generic fibre of $R^\univ$. If they were
to not form a closed set then $R_\tau$ would have quotients
corresponding to points which were not of type $\tau$, but which were
``close'' to being of type $\tau$ (more precisely, whose reductions modulo~$p^n$
were also reductions of type $\tau$ representations).
The paper~\cite{cdt} calls the points in the closure ``weakly
of type $\tau$'' and conjectures that being weakly of type $\tau$
is equivalent to being of type $\tau$. This conjecture
was proved not long afterwards for tame types by David Savitt
in~\cite{savitt}.

Now of course, one hopes that for certain types $\tau$,
the corresponding rings $R_\tau$ are small enough for the Taylor--Wiles
method to work (subject to the restriction that $p>2$
and that $\rhobar$ is absolutely irreducible even when restricted
to the absolute Galois group of $\Q(\sqrt{(-1)^{(p-1)/2}p})$, an assumption
needed to make the Taylor--Wiles machine work),
and big enough to capture some new elliptic curves.
Even though the definition of $R_\tau$ is in some sense a little
convoluted, one can still hope to write down a surjection
$W(k)[[t]]\to R_\tau$ in some cases (and thus control the tangent
space of $R_\tau$), for example by writing
down a deformation problem which is known to be
representable by a ring isomorphic to $W(k)[[t]]$, and showing
that it contains all the points of type $\tau$ (geometrically,
we are writing down a closed subset of $\Spec(R^\univ)$ with
sufficiently small tangent space,
checking it contains all the points of type $\tau$ and
concluding that it contains all of $\Spec(R_\tau)$). The problem
with such a strategy is that it requires a good understanding
of finite flat group schemes over the integers of the $p$-adic
field~$K$ corresponding to the kernel of~$\tau$. In 1998
the only fields for which enough was known were those extensions
$K$ of~$\Q_p$ which were tamely ramified. For such extensions,
some explicit calculations were done in~\cite{cdt} at the primes~3
and~5, where certain explicit $R_\tau$ were checked to have small
enough tangent space. There is a general modularity lifting
theorem announced in~\cite{cdt} but it includes, in the non-ordinary
case, an assumption the statement that $R_\tau$ is small enough for the method
to work, and this is difficult to check in practice, so the result
has limited applicability. However the authors
did manage to check this assumption in several explicit cases
when $p\in\{3,5\}$, and deduced

\begin{theorem}\label{cdtt} If $E/\Q$ is an elliptic curve which becomes semistable
at~3 over a tamely ramified extension of $\Q_3$, then~$E$ is modular.
\end{theorem}

This is the main theorem of~\cite{cdt} (see the second page
of loc.\ cit.). The proof is as follows: if $E[3]$ is irreducible
when restricted to the absolute Galois group of $\Q(\sqrt{-3})$
then they verify by an explicit calculation that either $E$ is
semistable at~3, or some appropriate $R_\tau$ is small enough,
and in either case this is enough. If $E[5]$ is irreducible
when restricted to the absolute Galois group of $\Q(\sqrt{-5})$
then $E[5]$ can be checked to be modular via the $3$--$5$ trick,
and $E$ can be proven modular as a consequence, although again
the argument relies on computing enough about an explicit $R_\tau$
to check that it is small enough. Finally Noam Elkies checked
for the authors that the number of $j$-invariants of elliptic curves
over~$\Q$ for which
neither assertion holds is finite and worked them out explicitly;
each $j$-invariant was individually checked to be modular.

At around the same time, Christophe Breuil had proven the breakthrough
theorem~\cite{breuil:annals}, giving a ``linear algebra'' description
of the category of finite flat group schemes over the integers
of an \emph{arbitrary} $p$-adic field. Armed with this, Conrad,
Diamond and Taylor knew that there was a chance that further
calculations of the sort done in~\cite{cdt} had a chance of
proving the full Taniyama-Shimura conjecture. The main problem
was that the rings $R_\tau$ were expected to be small enough
for quite a large class of tame types $\tau$, but were rarely
expected to be small enough if $\tau$ was wild. After much
study, Breuil, Conrad, Diamond and Taylor found an explicit
finite list of triples $(p,\rhobar,\tau)$ for which $R_\tau$ could
be proved to be small enough ($p$ was always~3 in this list,
and in one extreme case
they had to use a mild generalisation of a type called an ``extended
type'' in a case where $R_\tau$ was just too big; the extended type
cut it down enough), and this list and the 3--5 trick was enough
to prove

\begin{theorem}[Breuil, Conrad, Diamond, Taylor (2001)]\label{stw}
Any elliptic curve $E/\Q$ is modular.
\end{theorem}

\section{Kisin's modularity lifting theorems.}

In this section we briefly mention some important work of Kisin
that takes the ideas above much further.

As we have just explained, the Breuil--Conrad--Diamond--Taylor
strategy for proving a modularity lifting theorem was to write
down subtle local conditions at~$p$ which were representable
by a ring which was ``not too big'' (that is, its tangent
space is at most 1-dimensional). The main problem with
this approach was that the rings that this method needs to use
in cases where the representation is coming from a curve of large
conductor at~$p$ are (a) difficult to control, and (b) very rarely small enough
in practice. The authors of~\cite{bcdt} only just got away
with proving modularity of all elliptic curves because of some
coincidences specific to the prime~$3$, where the rings turned
out to be computable using Breuil's ideas, and just manageable enough
for the method to work. These calculations inspired conjectures
of Breuil and~M\'ezard~(\cite{bm}) relating an invariant of $R_\tau$
(the Hilbert--Samuel multiplicity of the mod~$p$ reduction
of this ring) to a representation-theoretic invariant (which
is much easier to compute).

Kisin in the breakthrough paper~\cite{kisin:annals} (note that this paper
was published in~2009 but the preprint had been available
since~2004) gave a revolutionary new way to approach the problem
of proving modularity lifting theorems. Kisin realised that rather
than doing the commutative algebra in the world of $\Z_p$-algebras,
one could instead just carry around the awkward rings $R_\tau$
introduced in~\cite{cdt}, and instead do all the dimension-counting
in the world of $R_\tau$-algebras (that is, count relative dimensions
instead). This insight turns out to seriously reduce
the amount of information one needs
about $R_\tau$; rather than it having to have a 1-dimensional tangent space,
it now basically only needs to be an integral domain of Krull
dimension~2. In fact one can get away with even less (which is good
because $R_\tau$ is not always an integral domain); one can even
argue using only an irreducible component of $\Spec(R_\tau[1/p])$, as long as
one can check that the deformations one is interested in live on
this component. 

There is one problem inherent in this method, as it stands: the
resulting modularity lifting theorems have a form containing
a condition which might be tough to verify in practice. For
example, they might say something like this: ``say $\rhobar$
is modular, coming from a modular form~$f$. Say $\rho$ lifts
$\rhobar$. Assume furthermore that $\rho_f$ and $\rho$ both
correspond to points on the same component of some $\Spec(R_\tau[1/p])$.
Then $\rho$ is modular.'' The problem here is that one now needs
either to be able to check which component various deformations
of $\rhobar$ are on, or to be able somehow to jump between components
(more precisely, one needs to prove theorems of the form ``if $\rhobar$
is modular coming from some modular form, then it is modular coming
from some modular form whose associated local Galois representation
lies on a given component of $\Spec(R_\tau[1/p])$''. 
Kisin managed to prove that certain $R_\tau$ only had one component,
and others had two components but that sometimes one could move
from one to the other, and as a result of these ``component-hopping''
tricks ended up proving the
following much cleaner theorem (\cite{kisin:annals}):

\begin{theorem}[Kisin]\label{kisinQ}
Let $p>2$ be a prime, let $\rho$ be a 2-dimensional
$p$-adic representation of $\GQ$ unramified outside a finite set of primes,
with reduction $\rhobar$, and assume that $\rhobar$ is modular
and $\rhobar|\Gal(\Qbar/K)$ is absolutely irreducible, where
$K=\Q(\sqrt{(-1)^{(p-1)/2}p})$. Assume furthermore that $\rho$
is potentially Barsotti--Tate and has determinant equal to
a finite order character times the cyclotomic character.
Then $\rho$ is modular.
\end{theorem}

Note that we do not make any assumption on the type of $\rho$; this
is why the theorem is so strong. This result gives another proof
of the modularity of all elliptic curves, because one can argue at~3
and~5 as in~\cite{cdt} and~\cite{bcdt} but is spared the hard computations
of $R_\tau$ in~\cite{bcdt}: the point is that Kisin's machine can often
deal with them even if their tangent space has dimension greater
than~one by doing the commutative algebra in this different and more
powerful way. These
arguments ultimately led to a proof of the Breuil--M\'ezard conjectures:
see for example Kisin's recent ICM talk.

The next part of the story in the case of 2-dimensional
representations of $\GQ$ would be the amazing work of Khare and
Wintenberger~(\cite{kw1}, \cite{kw2}), proving

\begin{theorem}[Khare--Wintenberger] Serre's conjecture (Conjecture~\ref{serreconj}) is true.
\end{theorem}
We have to stop somewhere however, so simply refer
the interested reader to the very readable
papers~\cite{kw1} and~\cite{khare:strat}
for an overview of the proof of this breakthrough result.

The work of Khare and Wintenberger means that nowadays we do not
have to rely on the Langlands--Tunnell
theorem to ``get us going'', and indeed we now get two more proofs
of Fermat's Last theorem
and of the modularity of all elliptic curves: firstly, given an elliptic
curve~$E$, we can just choose a random large prime, apply Khare--Wintenberger
to $E[p]$ and then apply the theorem of Kisin above. Secondly, given
an elliptic curve, we can apply the Khare--Wintenberger theorem to
$E[p]$ for all~$p$ at once, and then use known results about level
optimisation in Serre's conjecture to conclude again that~$E$
is modular. In particular we get a proof of FLT that avoids
non-Galois cubic base change. However it seems to the
author that things like cyclic base change and
the Jacquet--Langlands theorem will still be
essentially used in this proof,
and hence even now it seems that to understand a full proof
of FLT one still needs to understand both a huge amount of
algebraic geometry and algebra, and also a lot of hard analysis.

\section{Generalisations to totally real fields.}

So far we have restricted our discussion to modularity lifting
theorems that applied to representations of the absolute
Galois group of~$\Q$. It has long been realised that even if
one is mainly interested in these sorts of questions over $\Q$, it is
definitely worthwhile to prove as much as one can for a general
totally real field, because then one can use base change tricks
(the proofs of which are in~\cite{langlands:bc} and use
a lot of hard analysis) to get more information about the situation over $\Q$.
One of the first examples of this phenomenon, historically, was the result
in \S0.8 of~\cite{carayol:conductor}, where Carayol proves that the conductor
of a modular elliptic curve over $\Q$ was equal to the level of
the newform giving rise to the curve---even though this is a statement
about forms over $\Q$, the proof uses Hilbert modular forms over
totally real fields.

Some of what we have said above goes through to the totally real
setting. Let us summarise the current state of play. We fix a totally
real number field~$F$. The role of modular forms in the previous
sections is now played by Hilbert modular forms; to a Hilbert
modular eigenform there is an associated 2-dimensional $p$-adic
representation of the absolute Galois group of~$F$, and this representation
is totally odd, in the sense that the determinant of $\rho(c)$
is $-1$ for all complex conjugations~$c$ (there is more than one conjugacy class
of such things if $F\not=\Q$, corresponding to the embeddings $F\to\R$).
So formally the situation is quite similar to the case of $F=\Q$.
``Under the hood'' there are some subtle differences, because there
is more than one analogue of the theory of modular curves in this
setting (Hilbert modular varieties and Shimura curves, both of which
play a role, as do certain $0$-dimensional ``Shimura varieties''), but we
will not go any further into these issues; the point is that one
can \emph{formulate} the notion of modularity, and hence
of modularity lifting theorems in this setting.
But can one prove anything? One thing we certainly cannot prove,
at the time of writing, is

\begin{conjecture}[``Serre''] Any continuous totally odd irreducible
representation $\rhobar:\GF\to\GL_2(\Fpbar)$ is modular.
\end{conjecture}
Serre did not (as far as we know) formulate his conjecture in this generality,
but it has become part of the folklore and his name
seems now to be attached to it. We mention this conjecture because
of its importance in the theory. If one could prove this sort of
conjecture then modularity
of all elliptic curves over all totally real fields would
follow.

Just as in the case of $F=\Q$, the conjecture can be refined---for
example one can predict the weight and level of a modular form
that should give rise to $\rhobar$, as Serre did for $F=\Q$
in~\cite{serre:conj}. These refinements, in the case of the
weight, can be rather subtle: we refer the reader to~\cite{bdj}
when~$F$ is unramified at~$p$, and to work of Michael Schein
(for example~\cite{schein:ram})
and Gee (Conjecture~4.2.1 of~\cite{gee:prescribed}).

Once these refinements are made, one can ask two types of questions. The
first is of the form ``given $\rhobar$, is it modular?''. These sorts
of questions seem to be wide open for a general totally real field.
The second type is of the
form ``given $\rhobar$ which is assumed modular, is it modular
of the conjectured weight and level?''. This is the sort of question
which was answered by Ribet (the level: see~\cite{ken:100})
and Edixhoven (the weight: see~\cite{bas:weight}), following work of
many many others, for $F=\Q$. Much progress has also been made
in the general totally real case. For example see work
of Jarvis~\cite{jarvis} and Rajaei~\cite{rajaei} on the level,
and Gee~\cite{gee:wtopt} on the weight, so the situation
for Serre's conjecture on Hilbert modular forms now is becoming basically
the same as it was in the classical case before Khare--Wintenberger:
various forms of the conjecture are known to be equivalent, but all
are open.

Given that there is a notion of modularity, one can formulate
modularity lifting conjectures in this setting. But what can one
prove? The first serious results in this setting were produced by
Skinner and Wiles
in~\cite{sw1} and~\cite{sw2}, which applied in the setting of
\emph{ordinary} representations (that is, basically, to representations
$\rho$ which were upper triangular when restricted to a decomposition group for
each prime above~$p$). Here is an example of a modularity lifting
theorem that they prove (\cite{sw2} and correction in \cite{skinnerletter}):

\begin{theorem}[Skinner--Wiles, 2001] Suppose $p$ is an odd prime
and $F$ is a totally real field. Suppose $\rho:\GF\to\GL_2(\Qpbar)$
is continuous, irreducible, and unramified outside a finite set of places.
Suppose that $\det(\rho)$ is a finite order character times
some positive integer power of the cyclotomic character,
that $\rho|D_\p$ is upper triangular with an unramified quotient,
for all $\p|p$, and that the two characters on the diagonal
are distinct modulo $p$. Finally, suppose that $\rhobar|G_{F(\zeta_p)}$ is
absolutely irreducible, and that $\rhobar$ is modular,
coming from an \emph{ordinary} Hilbert modular form~$f$ of parallel
weight such that, for all $\p|p$, the unramified quotients of $\rho_f|D_\p$
and $\rho|D_\p$ are congruent mod~$p$. Then $\rho$ is modular.
\end{theorem}

Note that Skinner and Wiles show that if $\rhobar$ has an ordinary
modular lift, then many of its ordinary lifts are modular. Their
technique is rather more involved than the usual numerical criterion
argument---they make crucial use of deformations to characteristic~$p$
rings, and in fact do not show that the natural map $R\to T$ is an
isomorphism, using base change techniques to reduce to a case
where they can prove that it is a surjection with nilpotent kernel.

Kisin's work on the rings $R_\tau$ of the previous section all
generalised to the totally real setting, enabling Kisin to prove
some stronger modularity lifting theorems which were not confined
to the ordinary case. Kisin's original work required~$p$ to be
totally split in~$F$, but Gee proved something in the general
case. We state Gee's theorem below.

\begin{theorem}[Gee \cite{gee1},\cite{gee2}] Suppose $p>2$, $F$ is totally
real, and $\rho$ is a continuous potentially
Barsotti--Tate 2-dimensional $p$-adic Galois representation of
the absolute Galois group of~$F$, unramified outside
a finite set of primes, and with determinant equal to a finite order
character times the cyclotomic character.
Suppose that its reduction $\rhobar$ is modular,
coming from a Hilbert modular form~$f$, and suppose that for
all $v|p$, if $\rho$ is potentially ordinary at~$v$ then
so is $\rho_f$. Finally suppose $\rhobar$ is irreducible
when restricted to the absolute Galois group of $F(\zeta_p)$,
and if $p=5$ and the projective image of $\rhobar$ is isomorphic
to $\PGL_2(\F_5)$ then assume furthermore than $[F(\zeta_5):F]=4$.

Then $\rho$ is modular.
\end{theorem}

Note that, in contrast to Kisin's result Theorem~\ref{kisinQ},
in this generality ``component hopping'' is not as easy,
and the assumption in this theorem that if $\rho$ is potentially
ordinary then $\rho_f$ is too, are precisely assumptions ensuring
that $\rho$ and $\rho_f$ are giving points on the same components
of the relevant spaces $\Spec(R_\tau[1/p])$.

It is also worth remarking here that the Langlands--Tunnell theorem,
Theorem~\ref{lt}, is true for totally odd irreducible
representations of any totally real field to $\GL_2(\F_3)$,
so we can start to put together what we have to prove
some modularity theorems for elliptic curves. Note that Gee's
result above has the delicate assumption that not only
is $\rhobar$ modular, but it is modular coming from a Hilbert
modular form whose behaviour at primes dividing~$p$ is similar
to that of $\rho$. However Kisin's ``component hopping'' can be done
if~$p$ is totally split in~$F$, and Kisin can, using
basically the same methods, generalise his Theorem~\ref{kisinQ} to
the totally real case if $p$ is totally split, giving the following
powerful modularity result:

\begin{theorem}[Kisin~\cite{kisin:book}] Let $F$ be a totally
real field in which a prime $p>2$ is totally split, let
$\rho$ be a continuous irreducible 2-dimensional representation of $\GF$,
unramified outside a finite set of primes, and potentially
Barsotti--Tate at the primes above~$p$. Suppose that $\det(\rho)$
is a finite order character times the cyclotomic character,
that $\rhobar$ is modular coming from a Hilbert modular
form of parallel weight 2, and that $\rhobar|G_{F(\zeta_p)}$ is
absolutely irreducible. Then $\rho$ is modular.
\end{theorem}

As a consequence, if $F$ is a totally real field in which $p=3$
is totally split, and if $E/F$ is an elliptic curve with $E[3]|F_{F(\zeta_3)}$
absolutely irreducible, then $E$ is modular. 

Of course \emph{all} elliptic curves over $F$ are conjectured to
be modular, but this conjecture
still remains inaccessible. If one were to attempt to mimic
the strategy of proof in the case $F=\Q$ then one problem
would be that for a general totally real field, there may
be infinitely many elliptic curves with subgroups of order~15
defined over~$F$, and how can one deal with such curves?
There are infinitely many, so one cannot knock them off
one by one as Elkies did. Their mod~$3$ and mod~$5$ Galois
representations are globally reducible, and the best modularity
lifting theorems we have in this situation are in~\cite{sw1},
where various hypotheses on~$F$ are needed (for example $F/\Q$ has
to be abelian in Theorem~A of~\cite{sw1}). On the other hand,
because Serre's conjecture is still open for totally real
fields one cannot use the $p$-torsion for any prime $p\geq7$ either,
in general. It is not clear how to proceed in this situation!

\section{Potential modularity pre-Kisin and the $\p$--$\lambda$ trick.}

We have been daydreaming in the previous section about the possibilities
of proving that a general elliptic curve $E$ over a general totally
real field~$F$ is modular, and observing that we are not ready to prove
this result yet. Modularity is a wonderful thing to know for an elliptic
curve; for example, the Birch--Swinnerton-Dyer conjecture is a statement
about the behaviour of the $L$-function of an elliptic curve at the
point $s=1$, but the $L$-function of an elliptic curve is defined
by an infinite sum which converges for $\re(s)>3/2$, and it is only
a conjecture that this $L$-function has an analytic continuation
to the entire complex plane. One very natural way of analytically
continuing the $L$-function is to prove that the curve is modular,
because modular forms have nice analytic properties and the
analytic continuation of their $L$-functions is well-known.

For an elliptic curve over a general number field though, the $L$-function
is currently not known to have an analytic continuation, or even a
meromorphic continuation! However, perhaps surprisingly, it turns out
that the results above, plus one more good new idea due to Taylor,
enabled him to prove \emph{meromorphic} continuation for a huge class of
elliptic curves over totally real fields, and the ideas have now
been pushed sufficiently far to show that the $L$-function
of every elliptic curve over every totally real field can be
meromorphically continued. We want to say something about how this
all happened.

The starting point was Taylor's paper~\cite{taylor:fm}.
The basic idea behind this breakthrough paper is surprisingly easy
to explain! Recall first Lemma~\ref{35trick}, the $3$--$5$ trick.
We have an elliptic curve~$E/\Q$ with $E[5]$ irreducible,
and we want to prove that $E[5]$ is modular. We do this by writing down
a second elliptic curve~$A/\Q$ with $A[3]$ irreducible
and $A[5]\cong E[5]$. Then the trick, broadly, was that $A[3]$ is modular
by the Langlands--Tunnell theorem, so $A$ is modular by a modularity
lifting theorem, so $A[5]$ is modular, so $E[5]$ is modular. The proof
crucially uses the fact that the genus of the modular curve $X(5)$
is zero so clearly does not generalise to much higher numbers.

However, if we are allowed to be more flexible with our base field,
then this trick generalises very naturally and easily. Let us say
that we have an elliptic curve~$E$ over a totally real field~$F$,
and we want to prove that~$E$ is potentially modular (that is,
that~$E$ becomes modular over a finite extension field~$F'$ of~$F$,
also assumed totally real). Here is a
strategy. Say $p$ is a large prime such that $E[p]$ is irreducible.
Let us write down a random odd 2-dimensional
mod~$\ell$ Galois representation $\rho_\ell:\GF\to\GL_2(\Flbar)$
($\ell$ a prime, $\ell\not=p$) which is induced from a character;
because this representation is induced it is known to be modular.
Now let us consider the moduli space parametrising elliptic curves~$A$
equipped with
\begin{enumerate}
\item an isomorphism $A[p]\cong E[p]$
\item an isomorphism $A[\ell]\cong\rho_\ell$.
\end{enumerate}
This moduli problem will be represented by some
modular curve, whose connected components will be twists of $X(p\ell)$
and hence, if $p$ and $\ell$ are large, will typically have large genus.
However, such a curve may well still have lots of rational points, as long
as I am allowed to look for such things over an arbitrary finite extension
$F'$ of~$F$! So here is the plan: first, consider this moduli problem.
Second, find a point on this moduli space defined over~$F'$, for $F'$
some finite extension of~$F$. Next, ensure that $F'$ is totally
real, and that our modularity lifting theorems are robust
enough to apply in the two situations in which we will need them. More
precisely, we need one modularity lifting theorem of
the form ``$\rho_\ell$ is modular over $F'$ and hence $A/F'$ is modular'',
so $A[p]=E[p]$ is modular, and then another one which says
``$E[p]$ is modular over $F'$, and hence $E$ is modular''.

In~2001, our knowledge of modularity lifting theorems
was poorer than it is now (because, for example, this was
the era of~\cite{bcdt} and before Kisin's work on local deformation
problems), so this idea would not run as far as it naturally
wanted to. Let us sketch some of the issues that arise
here. Firstly, if our moduli space has no real points at all
for some embedding $F\to\R$ then we cannot find a point over
\emph{any} totally real extension. So we need to check our moduli problem
has real points. Secondly, the modularity lifting theorems available to
Taylor in the totally real case were those of Skinner and Wiles
so only applied in the ordinary setting (which is not much of
a restriction because the curve will be ordinary at infinitely many
places) but furthermore only applied in the ``distinguished''
setting (that is, the characters on the diagonal of the mod~$p$
local representation have to be distinct), so the completions of $F'$
at the primes above~$p$ had better not be too big, and similarly
the completions of~$F'$ at the primes above $\ell$ must not be
too big either. This results in more local conditions on $F'$, and
so we need to ensure two things: firstly that our moduli problems
have points defined over reasonably small extensions of $F_\p$
and $F_\lambda$ for $\p\mid p$ and $\lambda\mid\ell$, and secondly
that there is no local-global obstruction to the existence of a well-behaved
$F'$-point (that is: given that our moduli problem has points over certain
``small'' local fields, we need to ensure it has a point over a
totally real number field whose completions at the primes above $p$
and $\ell$ are equally ``small''). Fortunately, such a local-global theorem
was already a result of Moret-Bailly (see~\cite{mb}):

\begin{theorem}\label{mbthm} Let $K$ be a number field
and let~$S$ be a finite set of places of~$K$. Let~$X$
be a geometrically irreducible smooth quasi-projective
variety over~$K$. Let~$L$ run through the finite field
extensions of~$K$ in which all the primes of~$S$ split
completely. Then the union of $X(L)$, as $L$ runs through
these extensions, is Zariski-dense in~$X$.
\end{theorem}

Let us see how our plan looks so far. The idea now
is that given an elliptic curve $E/F$, we find a prime~$p$
such that the Skinner--Wiles theorem applies to $E[p]$
(so $p$ is an ordinary prime for which $E[p]$ locally
at~$p$ has two distinct characters on the diagonal)
and then write down a random odd prime~$\ell$
and an induced representation $\rho_\ell$
of $\GF$. Because Moret-Bailly's theorem needs $X$ geometrically
irreducible, let us ensure that $\rho_\ell$ has cyclotomic determinant,
and fix an alternating pairing on the underlying vector space
(to be thought of as a Weil pairing). Consider now the moduli space of elliptic
curves~$A$ equipped with isomorphisms $A[p]\cong E[p]$ and $A[l]\cong\rho_\ell$
both of which preserve the Weil pairing. If we can find points
on this curve defined over the completions of~$F$ at all primes
above~$p$, $\ell$ and~$\infty$, then we might hope to conclude.
But there are obstacles to finding such points. For example,
if $\lambda\mid\ell$ is a prime of~$F$
and $A[p]$ is unramified at~$\lambda$, and $b_\lambda$ is the number of points
on $A$ mod~$\lambda$, then one can read off $b_\lambda$ mod~$p$ from
$A[p]$, and one also has the Weil bounds on the integer $b_\lambda$,
and these two constraints on~$b_\lambda$ might not be simultaneously
satisfiable. Because of obstructions of this form, we do need to
slightly modify the strategy. Taylor
finds it easier to work not moduli spaces of elliptic curves over~$F$,
but with moduli spaces of so-called Hilbert--Blumenthal abelian
varieties over~$F$, that is, of higher-dimensional abelian varieties
(say $g$-dimensional) equipped with a certain kind of polarization
and endomorphisms by the integers in a second totally
real field~$M$ of degree~$g$ over~$\Q$. These more general abelian
varieties still give rise to 2-dimensional Galois representations,
and provide the extra flexibility necessary, at the expense
of making things technically more complicated---indeed, as is becoming clear,
this generalisation of the 3--5 trick is one of these ideas
where the great inspiration now has to be offset by the
large amount of perspiration that has to get the idea to
work. On the other hand, the idea is certainly not restricted
to Galois representations coming from elliptic curves, and in fact
applies to a large class of ordinary
Galois representations. Taylor managed to put everything together,
even in the ``pre-Kisin'' modularity world,
and managed in~\cite{taylor:fm} to prove that a wide class
of ordinary 2-dimensional Galois representations of $\GQ$
were potentially modular. Let us state his result here.

\begin{theorem}[Taylor] Let $p$ be an odd prime, let $\rho:\GQ\to\GL_2(\Qpbar)$
be a continuous odd irreducible representation unramified outside
a finite set of primes, and assume
$$\rho|D_p=\begin{pmatrix}\chi^n\psi_1&*\\0&\psi_2\end{pmatrix}$$
with $\chi$ the cyclotomic character, $n\geq1$, and $\psi_1$ and
$\psi_2$ two finitely ramified characters, such that the mod~$p$
reductions of $\chi^n\psi_1$ and $\psi_2$ are not equal on the
inertia subgroup of $D_p$. Then $\rho$ becomes modular
over some totally real number field.
\end{theorem}

As we said, the reason for the ordinarity assumption
is that the result was proved in~2000 before the more recent breakthroughs
in modularity lifting theorems. Taylor went on in~2001 in~\cite{taylor:fm2}
to prove an analogous theorem in the low weight crystalline case.

\section{Potential Modularity after Kisin.}

In this last section we put together Kisin's modularity lifting
theorem methods with Taylor's potential modularity ideas.
Together, the methods can be used to prove much stronger results such as
the following:

\begin{theorem} Let $E/F$ be an elliptic curve over a totally real
field. Then there is some totally real extension $F'/F$
such that $E/F'$ is modular.
\end{theorem}

In particular, the $L$-function of $E$ has meromorphic continuation
to the whole complex plane. It is difficult to give a precise attribution
to this theorem---the history is a little complicated. Perhaps
soon after Taylor saw Kisin's work on local deformation rings he
realised that this theorem was accessible, but perhaps he could also see
that the Sato--Tate conjecture was accessible, and turned
his attention to this problem instead. Whatever the history, it seems
that it was clear to the experts around 2006 and possibly even earlier
that the theorem above was accessible. The first published proof
that we are aware of is in the appendix by Wintenberger to~\cite{nekovar},
published in~2009. Much has happened recently in higher-dimensional
generalisations of modularity lifting theorems---so-called automorphy
lifting theorems, proving that various $n$-dimensional Galois
representations are automorphic or potentially automorphic, and
as a result one could also point to, for example,
Theorem~8.7
of~\cite{blght}, where a far stronger ($n$-dimensional) result is proven
from which the theorem follows.
Another place to read about the details of the proof of this potential
modularity theorem would be
the survey article of Snowden~\cite{snowden}, who sticks to the
2-dimensional situation and does a very good
job of explaining what is needed. Snowden fills in various
gaps in the literature in order to make his paper relatively self-contained
modulo some key ideas of Kisin; if the reader looks at Snowden's paper
then they will see that the crucial ideas are basically due to
Taylor and Kisin. One of the problems that needs to be solved
is how to do the ``component-hopping'': a general modularity lifting
theorem might look like ``if $\rho_1$ is modular and $\rho_2$ is congruent
to $\rho_1$, and $\rho_1$, $\rho_2$ give points on the same
components of certain local deformation spaces, then $\rho_2$ is modular''.
As a result one needs very fine control on the abelian variety
that is employed to do the $\p$--$\lambda$ trick; however this fine
control can be obtained by Moret-Bailly's theorem. The interested
reader should read these references, each of which give the details
of the argument.

\section{Some final remarks.}

This survey has to stop somewhere so we just thought we would mention
a few things that we have not touched on. In the 2-dimensional
case there is of course the work of Khare and Wintenberger, which
we have mentioned several times but not really touched on more
seriously. As mentioned already, Khare and Wintenberger have done a
good job of summarising their strategy in their papers, and the
interested reader should start there. Another major area which
we have left completely untouched is the higher-dimensional
$R=T$ theorems in the literature, concerned with modularity
of $n$-dimensional Galois representations. Here one uses
automorphic forms on certain unitary groups to construct Galois
representations, and new techniques are needed. The first
big result is~\cite{cht}, proving an $R=T$ result at minimal
level. After Kisin's rings are factored into the
equation and issues with components are resolved one can
prove much stronger theorems; the state of the art at
the time of writing seems to be the preprint~\cite{blggt},
which is again very clearly written and might serve as a good
introduction to the area. We apologise for not saying
more about these recent fabulous works.

\providecommand{\bysame}{\leavevmode\hbox to3em{\hrulefill}\thinspace}
\providecommand{\MR}{\relax\ifhmode\unskip\space\fi MR }
\providecommand{\MRhref}[2]{%
  \href{http://www.ams.org/mathscinet-getitem?mr=#1}{#2}
}
\providecommand{\href}[2]{#2}

\end{document}